\documentclass[11pt,fullpage]{article}

\usepackage{amsmath}
\usepackage{graphicx}
\usepackage{amsfonts}
\usepackage{amssymb}
\usepackage[normalem]{ulem}
\usepackage[usenames]{color}
\numberwithin{equation}{section}
\usepackage{natbib}
 \bibpunct[, ]{(}{)}{,}{a}{}{,}%

\newtheorem{prop}{Proposition}
\newtheorem{theorem}{Theorem}

\newtheorem{corollary}{Corollary}
\newenvironment{proof}[1][Proof]{\textbf{#1:} }{\ \rule{0.5em}{0.5em}}

\renewcommand{\baselinestretch}{1.2}

\newcommand{\e}{\mathbb{E}}

\newcommand{\var}{\mbox{var}}
\newcommand{\tvar}{\mbox{v}\widetilde{\mbox a}\mbox{r}}
\newcommand{\Pro}{\mathbb{P}}
\newcommand{\cov}{\mathbb{C}\mbox{ov}}
\newcommand{\cL}{\mathcal{L}}
\newcommand{\cE}{\mathcal{E}}
\newcommand{\bR}{\mathbb{R}}

\title{\LARGE On a Single Server Queue Fed by a Scheduled Traffic with Pareto Perturbations    \vspace{0.5cm}
}
\author{{\large   Victor F. Araman \hspace{0.5cm}  Hong Chen \hspace{0.5cm} Peter W. Glynn\hspace{0.5cm} Li Xia}

\thanks{The first author is with the  Olayan School of Business, American University of Beirut, Beirut, va03@aub.edu.lb. The second
author is with the Shanghai Advanced Institute of Finance,  Shanghai Jiaotong University, Shanghai 200030, China, hchen@saif.sjtu.edu.cn. The third author is with the Management Science and Engineering department at Stanford University, Stanford, CA, 74305, glynn@stanford.edu. The fourth author is with the Business School at Sun Yat-Sen University, Guangzhou 510275, China, xiali5@sysu.edu.cn}}

\begin{document}
    \maketitle \thispagestyle{empty}
\date{}
\renewcommand{\thefootnote}{\fnsymbol{footnote}}
\setcounter{footnote}{1}
\renewcommand{\baselinestretch}{1.2} {\small\normalsize}

\begin{abstract}    
A``scheduled" arrival process is one in which the $n^{th}$ arrival
is scheduled for time $n$, but instead occurs at $n+\xi_n$ , where the
$\xi_j$'s are iid. We describe here the behavior of a single server queue fed by such traffic in which
the processing times are deterministic. A particular focus is on perturbation with Pareto-like tails but with finite mean. We obtain tail approximations for the steady-state workload in both cases where the queue is critically loaded and under a heavy-traffic regime. A key to our approach is our analysis of the tail behavior of a sum of independent Bernoulli random variables with parameters of the form $p_n\sim c \,n^{-\alpha}$ as $n\rightarrow\infty$, for $c>0$ and $\alpha>1$. 

\end{abstract}

\maketitle

\section{Introduction}
In conventional queueing models, it is frequently assumed that the
exogenous arrivals to the system are described by a renewal
(counting) process. Specifically, the sequence $\chi =(\chi
_{n}:n\geq 1)$ of inter-arrival times of successive customers is
assumed to be a sequence of independent and identically distributed
(iid) non-negative random variables (rv's). More complex (arrival)
traffic models can be obtained by assuming that the $\chi _{j}$'$s$
are Markov-dependent or form a stationary time series.

While such traffic models are frequently appropriate, there are some
modelling settings in which one may seek alternatives. One such
setting is that in which arrivals are scheduled in advance; for example, an outpatient clinic. Patients are typically
scheduled to arrive at regular fifteen or twenty minute intervals. Of
course, some patients arrive early for their appointments, and
others arrive late, so that there is some random variation present.
A natural traffic model to adopt here is to assume that the $n^{th}$
patient is scheduled to arrive at the clinic at time $nh,$
but actually arrives at time $nh+\xi _{n},$ where $\xi =(\xi
_{n}:n\geq 0)$ is a stationary sequence of rv's. We call such an
arrival process a ``scheduled traffic model'', and we refer to the
    $\xi _{n}$'$s$ as the (random) \textit{perturbations} about the schedule.
Given our terminology, it therefore seems natural to use the
notation $S/M/1,$ $S/G/1,$ and so forth to refer to a single-server
queue in which the arrivals follow a scheduled traffic model, and in
which the processing times are exponential, generally distributed,
etc.

Scheduled traffic is described in \cite{CoxSmith1961} as a possible
arrival model to a queue. However, no analysis is offered for it.
\cite{Chenetal02} consider an application of the
$S/D/1$ queue to the air traffic control space in the vicinity of an
airport, and show that the $S/D/1$ queue is frequently stable even
when $\rho=1.$ More recently, \citet{AramanGlynn12} show that the arrival counting process corresponding to a scheduled traffic arrival process with infinite mean Pareto-like perturbations converges to a fractional Brownian motion with $H<1/2$. They also obtain a heavy traffic limit theorem for a single server queue fed by such traffic. Our primary goal here is the analysis of queues fed by scheduled traffic when the perturbations have finite mean.

It turns out that in the context of heavy traffic theory, a single server queue fed by a scheduled traffic with  i.i.d. service times, (i.e., S/G/1 queue) behaves exactly the same as a D/G/1 queue. In that setting, the stochasticity of the i.i.d. random service times dominates the randomness present in the scheduled arrival process. This has been already established in the case where the perturbations have infinite mean (see,  \citet{AramanGlynn12}). 
We reach here a similar conclusion for the case of finite-mean perturbations, the details of which are relegated to the last section. Moreover, in a separate work, we show that the tail asymptotics of the waiting times in an S/G/1 queue is again the same as that of a D/G/1 queue. In view of all this, and in order to highlight and expose the impact of scheduled traffic on a queueing system, our primary focus in this paper is on a single server queue with \textit{deterministic} service times, (i.e., an S/D/1 queue.)

A recent paper by \citet{Honnappa18} studies what they call transitory queueing and, in this context, introduces scheduled traffic, as one case of particular interest.  Specifically, they consider a queueing system that faces a finite number of scheduled arrivals during a given finite horizon and restrict themselves to perturbations that are uniformly distributed. They develop fluid and diffusion limits by scaling the number of scheduled arrivals while keeping the horizon and the perturbations untouched.  By doing so, the perturbations become increasingly large relatively to the duration of the scheduled interarrival times. This leads to an asymptotic regime that is very different than ours.  We also note the work on optimized appointment scheduling in outpatient care (e.g., \citet{Zacharias17} and \citet{Kemper14}). These works introduce also uncertainty to an initial traffic that is deterministically scheduled. \citet{Zacharias17} model the process of taking an appointment, as well as the resulting in-clinic queueing that is then created. They assume that only a  fraction of those scheduled will show up, but those showing up follow a renewal-like process. \citet{Kemper14} suggest a procedure for  optimal regularly scheduled appointments that will be feeding a D/G/1 like queue. They discuss how to adapt this procedure when arrivals are perturbed by an i.i.d. r.v. that they assume to be substantially smaller than a typical job duration.

In the next section, we present some properties of scheduled traffic. We establish a close connection between scheduled traffic and sums of Bernoulli random variables. In Section~\ref{sec: Sums of Bern} we study the logarithmic and exact tail asymptotics for sums of independent Bernoulli random variables with probabilities of the form $p_n\sim c \,n^{-\alpha}$ as $n\rightarrow\infty$, for $c>0$ and $\alpha>1$. The results we obtain there allow us to infer the tail asymptotics of the arrival counting process associated with scheduled traffic. We next analyze the behavior of a single server queue when fed by a scheduled traffic. Specifically, in Sections~\ref{sec: Critical Loading SD1} and \ref{Sec:5_Heavy_Traffic} we investigate the $S/D/1$ queue, and obtain limiting results for the workload, both when the queue is critically loaded and under a heavy traffic regime. Finally, in Section~\ref{sec: 6_SG1}, we discuss the $S/G/1$ queue with random service times, and argue that under a heavy traffic regime, such a queue behaves identically to the corresponding $D/G/1$ queue.

\section{Properties of Scheduled Traffic}\label{Sec: Preliminaries}
Let $(\xi_j:j\in\mathbb{Z})$ be an i.i.d. sequence of perturbations. We note that independence of the $\xi_j$'s seems plausible in many settings, given that perturbation $j$ is typically determined by decisions or preferences that are idiosyncratic to consumer $j$. Given the $\xi_j$'s, we define the random measure $\widetilde{N}$ via
$$\widetilde{N}(A)=\sum_jI(j+\xi_j+U\in A),$$
where $U$ is a uniform r.v. on $[0,1]$ independent of the $\xi_j$'s. It is easily argued that $\widetilde{N}$ is time-stationary, in the sense that $\widetilde{N}(\cdot+t)\overset{D}=\widetilde{N}(\cdot)$ for $t\in\mathbb{R}$ (where $\overset{D}=$ denotes equality in distribution.) We further define the counting process $N=\big(N(t):t\geq 0\big)$ via $N(t)=\widetilde{N}\big((0,t]\big);$ $N(t)$ counts the cumulative number of arrivals to the system in $(0,t]$. Our focus, in this section, is on the scheduled arrival process $N$.

We start by noting that regardless of whether $\xi_0$ has infinite mean or not, $N$ is a unit intensity counting process. Specifically,
\begin{equation*}
\begin{aligned}
\mathbb{E}N(t)&=\sum_j\int_0^1\mathbb{P}\big(j+x+\xi_j\in(0,t]\big)dx\\
&=\sum_j\int_0^1\mathbb{P}\big(j+x+\xi_0\in(0,t]\big)dx\\
&=\int_{-\infty}^\infty\mathbb{P}\big(r+\xi_0\in(0,t]\big)dr\\
&=\mathbb{E}\int_{-\infty}^\infty I\big(r\in(-\xi_0,t-\xi_0]\big)dr\\
&=t.
\end{aligned}
\end{equation*}

In fact, regardless of the tails of the $\xi_j$'s, the counting process $N$ has light tails. In particular, the moment generating function of $N(t)$ is always finite-valued. Specifically, the independence of the $\xi_j$'s ensures that for any $\theta$,

\begin{equation*}
\begin{aligned}
\log\left(\mathbb{E}\exp(\theta\,N(t))\right)&=\sum_j\log\left(\int_0^1\e\exp\Big(\theta\, I\big(j+x+\xi_j\in(0,t]\big)\Big)\,dx\right)\\
&=\sum_j\log\left(\int_j^{j+1}\e\exp\Big(\theta\, I\big(r+\xi_0\in(0,t]\big)\Big)\,dr\right)\\
&=\sum_j\log\left(1+(e^\theta\, -1)\int_j^{j+1}\,\mathbb{P}(r+\xi_0\in(0,t])\,dr\right)\\
&\leq (e^\theta\, -1)\sum_j\int_j^{j+1}\,\mathbb{P}(r+\xi_0\in(0,t])\,dr=(e^\theta-1)\,t.
\end{aligned}
\end{equation*}

In order to obtain insight into the dependence structure of $N$, we next study its covariance properties. Set $\Delta N(t)=N(t)-N(t-1)$ for $t\geq 1$, and recall that
\begin{equation}
\begin{aligned}
&\cov(\Delta N(1), \Delta N(t))\\
&=\mathbb{E}\,\cov\big((\Delta N(1),\Delta N(t))|U\big)+\cov\big(\mathbb{E}(\Delta N(1)|U),\mathbb{E}(\Delta N(t)|U)\big);&\label{Eq: Cov}\end{aligned}
\end{equation}
see p. 392 of \cite{Ross15}. Noting that
\begin{align*}
  \Delta N(t) &=\sum_j I\big(j+\xi_j+U\in(t-1,t]\big) \\
  &= \sum_j I\big(j-\lfloor t\rfloor+\xi_j+U\in(t-\lfloor t\rfloor-1,t-\lfloor t\rfloor]\big) \\
  &\overset{D}=\sum_k I\big(k+\xi_k+U\in(t-\lfloor t\rfloor-1,t-\lfloor t\rfloor]\big),
\end{align*}
it is evident that $\mathbb{E}(\Delta N(t)|U)$ depends on $t$ only through $t-\lfloor t\rfloor$. The second term in (\ref{Eq: Cov}) does not decay to zero as $t\rightarrow\infty$ and it reflects the correlation due to the common random placement of the time origin associated with $U$. The more informative term on the right-hand side of (\ref{Eq: Cov}) is $\cov\big((\Delta N(1),\Delta N(t))|U\big)$. Note that for $t\geq 2$,

\begin{equation}
\begin{aligned}
&\cov\big((\Delta N(1),\Delta N(t))|U\big)\\
&=\sum_{i,j}\mathbb{P}(i+\xi_i+U\in(0,1],j+\xi_j+U\in(t-1,t]|U)\\
&~~~~-\sum_{i, j}\mathbb{P}(i+\xi_i+U\in(0,1])\,\mathbb{P}(j+\xi_j+U\in(t-1,t]|U)\\
&=\sum_{i\neq j}\mathbb{P}(i+\xi_i+U\in(0,1],j+\xi_j+U\in(t-1,t]|U)\\
&~~~~-\sum_{i,j}\mathbb{P}(i+\xi_i+U\in(0,1])\,\mathbb{P}(j+\xi_j+U\in(t-1,t]|U)\\
&=-\sum_{i}\mathbb{P}(i+\xi_0+U\in(0,1]|U)\,\mathbb{P}(i+\xi_0+U\in(t-1,t]|U),
\end{aligned}
\label{Eq: Cov_1st_term}\end{equation}so the conditional covariance is always non-positive. This is intuitively reasonable, since scheduled traffic has the characteristic that if an abnormally large number of customers arrive in an interval, this reduces the number available to arrive in a subsequent interval. We can now use (\ref{Eq: Cov_1st_term}) to develop asymptotics for the conditional covariance.\\

\begin{prop}\label{Prop: Cov}~
  \begin{itemize}
    \item[i.)] Suppose that $\xi_0$ has a bounded density for which there exists positive constants $c_1,c_2,\alpha_1,\alpha_2$ such that
    \begin{align*}
      &f(x)  \sim c_1\,x^{-\alpha_1-1}, \\
      &f(-x)  \sim c_2\,x^{-\alpha_2-1}
    \end{align*}
    as $x\rightarrow\infty$. If $\alpha_1<\alpha_2,$ then $$\cov\big((\Delta N(1),\Delta N(n))|U\big)\sim -c_1\,n^{-\alpha_1-1}$$
    as $n\rightarrow \infty$, whereas if $\alpha_2<\alpha_1,$ then $$\cov\big((\Delta N(1),\Delta N(n))|U\big)\sim -c_2\,n^{-\alpha_2-1}$$
    as $n\rightarrow \infty$. If $\alpha_1=\alpha_2,$ then $$\cov\big((\Delta N(1),\Delta N(n))|U\big)\sim -(c_1+c_2)\,n^{-\alpha_1-1}$$
    as $n\rightarrow \infty$.
    \item[ii.)]Suppose that $\xi_0$ has a bounded density $f$ for which there exists positive constants $d_1,d_2,\beta_1,\beta_2$ such that
    \begin{align*}
      &f(x)  \sim d_1\,e^{-\beta_1\,x}, \\
      & f(-x)  \sim d_2\,e^{-\beta_2\,x}
    \end{align*}
    as $x\rightarrow\infty$. If $\beta_1<\beta_2,$ then $$\cov\big((\Delta N(1),\Delta N(n))|U\big)\sim - e^{-\beta_1\,n}\,\frac{d_1}{\beta_1}\,(e^{\beta_1}-1)\sum_je^{-\beta_1 \,(j-U)}\,\mathbb{P}(\xi_0+U\in(j-1,j]|U)$$
    as $n\rightarrow \infty$, whereas if $\beta_2<\beta_1,$ then $$\cov\big((\Delta N(1),\Delta N(n))|U\big)\sim -e^{-\beta_2\,n}\,\frac{d_2}{\beta_2}\,(1-e^{-\beta_2})\sum_je^{-\beta_2 \,(j+U)}\,\mathbb{P}(\xi_0+U\in(j-1,j]|U)$$
    as $n\rightarrow \infty$. If $\beta_1=\beta_2,$ then $$\cov\big((\Delta N(1),\Delta N(n))|U\big)\sim -n\,e^{-\beta_1\,(n+1)}\frac{d_1\,d_2}{\beta_1^2}\,(1-e^{-\beta_1})^2$$
    as $n\rightarrow \infty$.
  \end{itemize}
\end{prop}
\hfill \break
\begin{proof}
  According to (\ref{Eq: Cov_1st_term}), the conditional covariance is given by
\begin{equation}
\begin{aligned}\nonumber
&-\sum_{j}\mathbb{P}(\xi_0+U\in(j-1,j]|U)\,\mathbb{P}(\xi_0+U\in(n+j-1,n+j]|U)\\
=&-\sum_{j>-n/2}\mathbb{P}(\xi_0+U\in(j-1,j]|U)\,\mathbb{P}(\xi_0+U\in(n+j-1,n+j]|U)\\
&-\sum_{k\leq n/2}\mathbb{P}(\xi_0+U\in(k-n-1,k-n]|U)\,\mathbb{P}(\xi_0+U\in(k-1,k]|U)
\end{aligned}
\end{equation}
Given our bounded density assumption, the Bounded Convergence Theorem implies that  $$n^{\alpha_1+1}\,\mathbb{P}(\xi_0+U\in(n+j-1,n+j])\rightarrow c_1$$
as $n\rightarrow\infty$, and $\big(n^{\alpha_1+1}\,\mathbb{P}(\xi_0+U\in(n+j-1,n+j]):j>-n/2\big)$ is uniformly bounded. Another application of the Bounded Convergence Theorem therefore implies that,
$$n^{\alpha_1+1}\sum_{j>-n/2}\mathbb{P}(\xi_0+U\in(j-1,j])\,\mathbb{P}(\xi_0+U\in(n+j-1,n+j]|U)\rightarrow c_1$$
as $n\rightarrow\infty$. Similarly,
$$n^{\alpha_2+1}\sum_{k\leq n/2}\mathbb{P}(\xi_0+U\in(k-n-1,k-n]|U)\,\mathbb{P}(\xi_0+U\in(k-1,k]|U)\rightarrow c_2$$
as $n\rightarrow\infty$, proving part i.).

For part ii.), suppose first that $\beta_1<\beta_2$ and note that $$\sum_je^{-\beta_1 j}\,\mathbb{P}(\xi_0+U\in(j-1,j]|U)<\infty.$$
Furthermore, our assumption on $f$ guarantees that
\begin{align*}
e^{\beta_1 n}\,\mathbb{P}(\xi_0+U\in(n+j-1,n+j]|U) \rightarrow \frac{d_1}{\beta_1}\,(e^{\beta_1}-1)\,e^{-\beta_1 (j-U)}
\end{align*}
as $n\rightarrow\infty$, and $\big(e^{\beta_1 j}\,\mathbb{P}(\xi_0+U\in(j-1,j]|U):j\geq 0\big)$ is uniformly bounded. Applying the Bounded Convergence Theorem, we conclude that

\begin{align*}
e^{\beta_1\,n}\,&\sum_{j>-n}\mathbb{P}(\xi_0+U\in(j-1,j]|U)\,\mathbb{P}(\xi_0+U\in(n+j-1,n+j]|U)\\
=&\sum_{j>-n}\mathbb{P}(\xi_0+U\in(j-1,j]|U)\,e^{-\beta_1 j}\cdot e^{\beta_1\,(n+j)}\,\mathbb{P}(\xi_0+U\in(n+j-1,n+j]|U)\\
&\rightarrow \frac{d_1}{\beta_1}\,(e^{\beta_1}-1)\,e^{\beta_1 U}\sum_j e^{-\beta_1 j}\,\mathbb{P}(\xi_0+U\in(j-1,j]|U)
\end{align*}
as $n\rightarrow\infty$. Similarly,
\begin{equation}
\begin{aligned}
e^{\beta_2\,n}\,&\sum_{k<0}\mathbb{P}(\xi_0+U\in(k-n-1,k-n]|U)\,\mathbb{P}(\xi_0+U\in(k-1,k]|U)\\
&\rightarrow \frac{d_2}{\beta_2}\,(1-e^{-\beta_2})\,e^{-\beta_2 U}\sum_{k<0}\mathbb{P}(\xi_0+U\in(k-1,k]|U)\,e^{\beta_2 \,k}\,
\end{aligned}
\label{Eq Cov_Term_Notes2.3}\end{equation}
as $n\rightarrow\infty,$ thereby establishing that the conditional covariance satisfies

\begin{align*}
e^{\beta_1 n}\,&\sum_{j}\mathbb{P}(\xi_0+U\in(j-1,j]|U)\,\mathbb{P}(\xi_0+U\in(n+j-1,n+j]|U)\\
&\rightarrow \frac{d_1}{\beta_1}\,(e{^\beta_1}-1)\,\sum_je^{-\beta_1 (j-U)}\,\mathbb{P}(\xi_0+U\in(j-1,j]|U)
\end{align*}
as $n\rightarrow\infty.$ The case where $\beta_2<\beta_1$ can be handled identically.

To handle the case where $\beta_1=\beta_2$, we write the conditional covariance as
\begin{align*}
&-\sum_{j\geq0}\mathbb{P}(\xi_0+U\in(j-1,j]|U)\,\mathbb{P}(\xi_0+U\in(n+j-1,n+j]|U)\\
&-\sum_{-n\leq j<0}\mathbb{P}(\xi_0+U\in(j-1,j]|U)\,\mathbb{P}(\xi_0+U\in(n+j-1,n+j]|U)\\
&-\sum_{k<0}\mathbb{P}(\xi_0+U\in(k-n-1,k-n]|U)\,\mathbb{P}(\xi_0+U\in(k-1,k]|U).
\end{align*}
Relation (\ref{Eq Cov_Term_Notes2.3}) shows that the third term is of order $O(e^{-\beta_2\,n})$ as $n\rightarrow\infty;$ a similar argument proves that the first term is of order $O(e^{-\beta_1\,n})$. To handle the second term, we write it as
\begin{equation}
\begin{aligned}
&-\sum_{-n\leq j<-n/2}\mathbb{P}(\xi_0+U\in(j-1,j]|U)\,\mathbb{P}(\xi_0+U\in(n+j-1,n+j]|U)\\
&-\sum_{-n/2\leq j<0}\mathbb{P}(\xi_0+U\in(j-1,j]|U)\,\mathbb{P}(\xi_0+U\in(n+j-1,n+j]|U)\\
=&-\sum_{0\leq k<n/2}\mathbb{P}(\xi_0+U\in(k-n-1,k-n]|U)\,\mathbb{P}(\xi_0+U\in(k-1,k]|U)\\
&-\sum_{-n/2\leq j<0}\mathbb{P}(\xi_0+U\in(j-1,j]|U)\,\mathbb{P}(\xi_0+U\in(n+j-1,n+j]|U).
\end{aligned}\label{Eq Cov_Term_Notes2.4}
\end{equation}
But the second term above equals
$$-\sum_{-n/2\leq j<0}\mathbb{P}(\xi_0+U\in(j-1,j]|U)\,\frac{d_1}{\beta_1}\,(e^{\beta_1}-1)\,e^{-\beta_1(n+j-U)}\,\big(1+o(1)\big),$$
as $n\rightarrow\infty$, where the $o(1)$ term is uniform in $-n/2\leq j<0.$ So this sum equals
$$-\big(1+o(1)\big)\,e^{-\beta_1\,n}\,\sum_{-n/2\leq j<0}\mathbb{P}(\xi_0+U\in(j-1,j]|U)\,\frac{d_1}{\beta_1}\,(e^{\beta_1}-1)\,e^{-\beta_1 (j-U)}.$$
Since $\beta_1=\beta_2,$ $\mathbb{P}\big(\xi_0+U\in(j-1,j]|U\big)\,e^{\beta_1 j}\rightarrow \frac{d_2}{\beta_1}\,(1-e^{-\beta_1})\,e^{-\beta_1 U}$ as $j\rightarrow-\infty.$ Consequently, the second term is asymptotic to $-e^{-\beta_1\,(n+1)}(n/2)\,\frac{d_1\,d_2}{\beta_1^2}\,(1-e^{-\beta_1})^2$ as $n\rightarrow\infty$. A similar analysis works for the first term in (\ref{Eq Cov_Term_Notes2.4}), proving part ii.) for $\beta_1=\beta_2.$
\end{proof}\\

As a consequence of Proposition~\ref{Prop: Cov} \textit{i.)}, we see that if $\min\{\alpha_1,\alpha_2\}\leq 1,$ the conditional autocorrelations are non-summable, indicating long-range dependence. This is the parameter range covered by Araman and Glynn (2012), in which it was established that $N(\cdot)$ satisfies a functional limit theorem with fractional Brownian motion having $H<1/2$ as a limit.

We turn next to a key representation for $N$ that holds only when $\mathbb{E}|\xi_0|<\infty$. In preparation for stating this result, let
\begin{align*}
\cE(t)&=\sum_{i+U>t}I(i+\xi_i+U\leq t)\\
\cL(t)&=\sum_{i+U\leq t}I(i+\xi_i+U> t).
\end{align*}
The r.v. $\cE(t)$ represents the total number of early customers at time $t$, who have arrived earlier than scheduled, while $\mathcal{L}(t)$ is the total number of late customers that will arrive after $t$ but were scheduled to arrive before $t$. The Borel-Cantelli lemma makes clear that $\cE(t)$ is finite-valued a.s. if and only if $\e \xi_0^- \overset{\Delta}=\e\max(-\xi_0,0)<\infty$ while $\cL(t)$ is finite-valued a.s. if and only if  $\e \xi_0^+ \overset{\Delta}=\e\max(\xi_0,0)<\infty$. Furthermore, $\big((\cE(t),\cL(t)):t\in\bR\big)$ is a time-stationary process,  where for every $t$, $\cE(t)$ and $\cL(t)$ are independent random variables.\\
\begin{prop}\label{Prop: N(t)-t}Suppose that $\mathbb{E}|\xi_0|<\infty$. Then, for $t \geq 0,$
\begin{equation}\label{Eq: N(t)-t}
  N(t)-t=\big(\sum_{i+U\in(0,t]}1\big)-t+\big(\cE(t)-\cL(t)\big)-\big(\cE(0)-\cL(0)\big).
\end{equation}
\end{prop}
\hfill \break
\begin{proof}
Observe that
\begin{equation}
\begin{aligned}
N(t)-t=&\sum_{i+U\in(0,t]}I(i+\xi_i+U\in(0,t])+\sum_{i+U>t}I(i+\xi_i+U\in(0,t])\\
&+\sum_{i+U\leq 0}I(i+\xi_i+U\in(0,t])-t\\
=&\sum_{i+U\in(0,t]}\big(1-I(i+\xi_i+U>t)-I(i+\xi_i+U\leq 0)\big)-t\\
&+\sum_{i+U>t}\big(I(i+\xi_i+U\leq t)-I(i+\xi_i+U\leq 0)\big)\\
&+\sum_{i+U\leq 0}I(i+\xi_i+U>0)-I(i+\xi_i+U>t).
\end{aligned}
\end{equation}
We now combine the first indicator sum with the sixth (to obtain $-\cL(t))$, and the second indicator sum with the fourth (to obtain $-\cE(0))$, thereby proving the result.
\end{proof}\\


We can now prove that $N(t)-t$ converges weakly as $t\rightarrow\infty,$ when we let $t\rightarrow\infty$ in such a way that $t-\lfloor t\rfloor$ is constant.

\begin{theorem}\label{Th: N(n)-n}
  Suppose that $\mathbb{E}|\xi_0|<\infty$, and fix $s\in[0,1).$ Then,
  \begin{equation*}
    N(n+s)-(n+s)\Rightarrow-s+I(U\leq s)+\big(\cE'(s)-\cL'(s)\big)-\big(\cE(0)-\cL(0)\big)
  \end{equation*}
  as $n\rightarrow\infty,$ where $\cE'(s),\cL'(s),\cE(0),\cL(0)$ are independent of one another given $U$, and $\cE'(s)\overset{D}=\cE(0)$, $\cL'(s)\overset{D}=\cL(0)$.
\label{Th: N(n)-n}\end{theorem}
\hfill \break
\begin{proof}
Recall that $$N(n+s)-(n+s)= \sum_{i+U\in(0,n+s]}1-(n+s)+(\cE(n+s)-\cE(0))-(\cL(n+s)-\cL(0)).$$
We start by observing that $$\sum_{i+U\in(0,n+s]}1-(n+s)=-s+I(U\leq s)$$ for $n\in\mathbb{Z}_+, s\in[0,1).$ Furthermore, if $k_n$ is an integer-valued sequence such that $k_n/n\rightarrow\nu\in(0,1)$ as $n\rightarrow\infty,$ we can write
\begin{align*}
&\big(\cE(n+s)-\cE(0),\cL(n+s)-\cL(0)\big)\\&=\Big(\sum_{i+U>n+s}I(i+U+\xi_i\in(0,n+s])-\sum_{i+U\in(0,n+s]}I(i+U+\xi_i\leq 0),\\
&~~~\sum_{i+U\in(0, n+s]}I(i+\xi_i+U> n+s)-\sum_{i+U\leq 0}I(i+U+\xi_i\in(0,n+s])
 \Big)\\
&=\Big(\sum_{i+U>n+s}I(i+U+\xi_i\in(0,n+s])-\sum_{i+U\in(0,k_n]}I(i+U+\xi_i\leq 0)-\sum_{i+U\in(k_n,n+s]}I(i+U+\xi_i\leq 0),\\
&~~~\sum_{i+U\in(0, k_n]}I(i+\xi_i+U> n+s)+\sum_{i+U\in(k_n, k_n+s]}I(i+\xi_i+U> n+s)-\sum_{i+U\leq 0}I(i+U+\xi_i\in (0,n+s])\Big)\\
\overset{D}=~& \big(\cE''(n+s)-\hat\cE_n-\hat\cE_n'',\hat\cL_n''(n+s)+\cL_n''(n+s)-\hat\cL_n\big).
\end{align*}
Note that because $\mathbb{E}\xi_0^+<\infty,$
\begin{align*}
  \mathbb{E}[\cL''_n(n+s)|U]& =\sum_{i+U\in(0, k_n]}\mathbb{P}(i+\xi_i+U>n+ s|U)\\
  & \leq \sum_{j+U\leq 0}\mathbb{P}(\xi_0>n-k_n-j-1)\rightarrow 0
\end{align*}
 as $n\rightarrow\infty$, proving that $\hat\cL''_n(n+s)\Rightarrow0$ as  $n\rightarrow\infty$. Similarly, the fact that $\mathbb{E}\xi_0^-<\infty$ implies that $\hat\cE_n''\Rightarrow 0$ as $n\rightarrow\infty.$ Finally, the four random variables $\big(\cE''(n+s),\hat\cE_n,\cL''_n(n+s),\hat\cL_n\big)$ all involve sums over subsets in $i$ that are disjoint, so they are conditionally independent of one another, given $U$. Furthermore,  $\cL''_n(n+s)\Rightarrow\cL'(s)$ and $\cE_n''(s)\Rightarrow \cE'(s)$ while, $\hat\cE_n\Rightarrow\cE(0)$ and $\hat\cL_n\Rightarrow\cL(0)$ as $n\rightarrow\infty$, proving the theorem. \end{proof}

Note that we must restrict convergence to sequences of the form $t_n=n+s$ with $n\rightarrow\infty$. In particular, weak convergence does not hold when $t\rightarrow\infty$ without any restrictions. To see this, consider the case in which $\xi_0=0~ a.s.$ Then,
 \begin{align*}
   N(t)-t &=\sum_{i+U\in(0,t]}1-t \\
   & =-(t-\lfloor t \rfloor)+I(U\leq t-\lfloor t \rfloor),
 \end{align*}
 and observe that the distribution of the right-hand side depends on  $t-\lfloor t \rfloor$, regardless of the magnitude of $t$.

 Theorem~\ref{Th: N(n)-n} shows that $N(t)-t$ is stochastically bounded in $t$. This is in sharp contrast to the case in which (for example) $N$ is a unit rate renewal counting process with finite-variance inter-arrival times, in which event $t^{-1/2}\big(N(t)-t\big)$ converges weakly to a normal r.v. (see Ross (1996)), so that $N(t)-t$ exhibits stochastic fluctuations of order $t^{1/2}.$

\section{Tail Asymptotics for Sums of Bernoulli Random Variables}\label{sec: Sums of Bern}
The analysis of Section~\ref{Sec: Preliminaries} establishes that $N,\cE$ and $\cL$ all can be clearly represented as sums of independent Bernoulli r.v.'s. As we will see in the next section, the tail behavior of these r.v.'s significantly affects the queueing dynamics of systems that are fed by scheduled traffic. In addition, Bernoulli sums arise in many other applications settings (e.g. credit risk). As a consequence, this section is focused on tail behavior for such Bernoulli sums.

Let $(I_j:j\in\mathbb{Z})$ be a family of independent r.v.'s, in which $p_j=\Pro(I_j=1)=1-\Pro(I_j=0).$
\begin{theorem}
Suppose that there exist constants $c>0$ and $\alpha>1$ for which $p_n\sim c\,n^{-\alpha}$ as $n\rightarrow\infty.$ If $Z=\sum_{j\geq 0}I_j,$ then $$\frac{1}{z\,\log z}\,\log\Pro(Z>z)\rightarrow-\alpha$$ as $z\rightarrow\infty.$\label{Th: Bernoulli Tail}
\end{theorem}
\hfill \break
\begin{proof}
  We shall employ an argument similar to that commonly used in the theory of large deviations; see, for example, p. 44 in \cite{Dembo}. (Note, however, that the asymptotic setting described by Theorem~\ref{Th: Bernoulli Tail} is not covered by traditional large deviations.) We start by observing that
  $$\psi(\theta)\overset{\Delta}=\log\e\exp\big(\theta\,Z\big)=\sum_{j\geq0}\log\big(p_j(e^\theta-1)+1\big)$$
  (where the sum converges absolutely since $\alpha>1).$ Choose $\theta=\theta(z)$ such that $e^{\theta(z)}=rz^\alpha$ (where $r>0$), and note that for $\epsilon>0$, $\theta>0,$ and $z$ sufficiently large,
\begin{align}
\psi\big(\theta(z)\big)=&\sum_{0\leq j\leq \lfloor \epsilon z\rfloor}\log\big(p_j(e^{\theta(z)}-1)+1\big)\nonumber\\
&+\sum_{j>\lfloor \epsilon z\rfloor}\log\big(p_j(e^{\theta(z)}-1)+1\big)\nonumber\\
\leq&\sum_{0\leq j\leq \lfloor \epsilon z\rfloor}\log\big(e^{\theta(z)}+1\big)\nonumber\\
&+\sum_{j>\lfloor \epsilon z\rfloor}\log\big((1+\epsilon)c\,j^{-\alpha}e^{\theta(z)}+1\big)\nonumber\\
&\leq  (\lfloor \epsilon z\rfloor+1)\,\log (1+rz^\alpha)+\sum_{j>\lfloor \epsilon z\rfloor}\log\big((1+\epsilon)\,rc\,(j/z)^{-\alpha}+1\big).\label{Eq: Bernoulli 3.1}
\end{align}
Note that the second term in (\ref{Eq: Bernoulli 3.1}), when multiplied by $1/z$, is a Riemann sum approximation, and hence
\begin{align*}
\frac{1}{z}\sum_{j>\lfloor \epsilon z\rfloor}\log&\big((1+\epsilon)\,rc\,(j/z)^{-\alpha}+1\big)\\
&\rightarrow\int_\epsilon^\infty\log\big((1+\epsilon)\,rc\,x^{-\alpha}+1\big)dx
\end{align*}
as $z\rightarrow\infty.$ (Specifically, the function $\log\big((1+\epsilon)\,rc\,x^{-\alpha}+1\big)$ is directly Riemann integrable (see \cite{Asmussen03}), so the Riemann approximation over $[\epsilon,\infty)$ converges.) It follows that
$$\overline{\lim}_{z\rightarrow\infty}\,\frac{1}{z\log z}\psi\big(\theta(z)\big)\leq \epsilon\,\alpha.$$

Markov's inequality guarantees that $$\Pro(Z>z)\leq \exp\big(-\theta(z)z+\psi(\theta(z))\big),$$ and hence
$$\overline{\lim}_{z\rightarrow\infty}\,\frac{1}{z\log z}\log\mathbb{P}(Z>z)\leq-\alpha\,(1-\epsilon).$$
Since $\epsilon>0$ can be chosen to be arbitrarily small, we conclude that
\begin{equation}\overline{\lim}_{z\rightarrow\infty}\,\frac{1}{z\log z}\log\mathbb{P}(Z>z)\leq-\alpha.\end{equation}

To obtain the lower bound needed for Theorem~\ref{Th: Bernoulli Tail}, we apply a change-of-measure argument. For $z>0$, put $$\widetilde{\Pro}_z(\cdot)=\e I(\cdot)\exp\big(\theta(z)\,Z-\psi(\theta(z))\big),$$ and let $\widetilde{\e}_z(\cdot)$ be the associated expectation operator. Then, \begin{equation}\Pro(Z>z)=\widetilde{\e}_zI(Z>z)\exp\big(-\theta(z)Z+\psi(\theta(z))\big).\label{Eq: Change_measure}\end{equation} Of course,
\begin{equation}\label{eq: Etilde_Z}\widetilde{\e}_zZ=\psi'(\theta(z))=\sum_{j\geq 0}\frac{p_je^{\theta(z)}}{p_j(e^{\theta(z)}-1)+1}.\end{equation} So,
$$\frac{1}{z}\,\widetilde{\e}_zZ=\sum_{j\geq 0}\frac{p_j\,rz^\alpha}{p_j(rz^\alpha-1)+1}\cdot \frac{1}{z}.$$ Since $p_j\,z^\alpha\sim c\,(j/z)^{-\alpha}$ as $j\rightarrow\infty,$ a simple adaptation of the earlier Riemann sum approximation argument proves that
$$\frac{1}{z}\,\widetilde{\e}_zZ\rightarrow\int_0^\infty\frac{c\,r}{c\,r+x^\alpha}dx$$ as $z\rightarrow\infty. $ Similarly,
\begin{align*}
\frac{1}{z}\,\tvar_zZ=\frac{1}{z}\psi''(\theta(z))&=\sum_{j\geq 0}\frac{p_j\,(1-p_j)\,e^{\theta(z)}}{(p_j\,(e^{\theta(z)}-1)+1)^2}\cdot\frac{1}{z}\\
&\rightarrow\int_0^\infty\frac{c\,r\,x^\alpha}{(c\,r+x^\alpha)^2}\,dx
\end{align*}
as $z\rightarrow\infty.$ For $\epsilon>0,$ we now select $r>0$ (uniquely) so that $$\int_0^\infty\frac{c\,r}{c\,r+x^\alpha}\,dx=1+\epsilon.$$ Observe that because $\psi(\theta(z))>0,$
\begin{align}
\Pro(Z>z)&=\widetilde{\e}_zI(Z>z)\exp\big(-\theta(z)\,Z+\psi(\theta(z))\big)\nonumber\\
&\geq \exp\big(-\theta(z)\,z+\psi(\theta(z))\big)\,\widetilde{\Pro}_z(Z>z)\nonumber\\
&\geq \exp\big(-\theta(z)\,z\big)\,\widetilde{\Pro}_z(Z\geq z).\label{Eq: Tail Bern Notes 3.3}
\end{align}
 But for $z$ large enough, we have that
\begin{align*}
\widetilde{\Pro}_z(Z>z)&= \widetilde{\Pro}_z(Z-\widetilde{\e}_zZ>z-\widetilde{\e}_zZ)\\
&\geq \widetilde{\Pro}_z(Z-\widetilde{\e}_zZ>z-(1+\varepsilon-\varepsilon/2)\,z)\\
&\geq \widetilde{\Pro}_z(Z-\widetilde{\e}_zZ>-\varepsilon\,z/2)\\
&\geq 1- \widetilde{\Pro}_z(|Z-\widetilde{\e}_zZ|>\epsilon\,z/2)\\
&\geq 1- 4 \frac{\tvar_zZ}{z^2\epsilon^2}\rightarrow 1
\end{align*}
as $z\rightarrow\infty,$ where the last inequality is an application of Chebyshev's inequality. Hence, (\ref{Eq: Tail Bern Notes 3.3}) implies that $$\underline{\lim}_{z\rightarrow\infty}\,\frac{1}{z\,\log z}\,\log\Pro(Z>z)\geq -\alpha,$$ proving the theorem.
\end{proof}\\

We now turn to the tail of $Z$ when $Z$ is the difference of two independent Bernoulli sums, $\sum_{j\geq0}I_j$ and $\sum_{j<0}{\tilde I}_j$.\\
\begin{corollary}\label{Cor: LD_Diff_Bern}
Suppose that $\mathbb{E}\sum_{j<0}{\tilde I}_j<\infty$ and that there exists $c>0$ and $\alpha>1$ for which $\mathbb{E}I_n\sim c\,n^{-\alpha}$ as $n\rightarrow+\infty$. If $Z=\sum_{j\geq 0}I_j -\sum_{j< 0}{\tilde I}_j$, then $$\frac{1}{z\log z}\log\Pro(Z>z)\rightarrow-\alpha$$ as $z\rightarrow\infty.$
\end{corollary}
\hfill \break
\begin{proof}
  We note that $\Pro(Z>z)\leq\Pro(\sum_{j\geq0}I_j>z)$, and apply Theorem~\ref{Th: Bernoulli Tail} to conclude that $$\overline{\lim}_{z\rightarrow\infty}\,\frac{1}{z\log z}\log\Pro(Z>z)\leq -\alpha.$$ For the lower bound, observe that the independence yields
  \begin{align*}
    \Pro(Z>z) & \geq\Pro\Big(\sum_{j\geq0}I_j>z+d,\sum_{j<0}I_j\leq d\Big) \\
     & =\Pro\Big(\sum_{j\geq0}I_j>z+d\Big)\,\Pro\Big(\sum_{j<0}I_j\leq d\Big).
  \end{align*}
  Hence, we apply Theorem~\ref{Th: Bernoulli Tail} to conclude that
  $$\underline{\lim}_{z\rightarrow\infty}\,\frac{1}{z\log z}\log\Pro(Z>z)\geq \underline{\lim}_{z\rightarrow\infty}\,\frac{1}{z\log z}\log \Pro\Big(\sum_{j\geq 0}I_j>z+d\Big)= -\alpha,$$ proving the result.
\end{proof}\\

We can immediately apply Theorem~\ref{Th: Bernoulli Tail} and its corollary to the tail asymptotics of $\cE,\cL$ and $N(t)$.\\
\begin{theorem}\label{th:Tail_N(t)}~
  \begin{itemize}
  \item[i.)] Suppose that $\xi_0$ is such that $\Pro(\xi_0>x)\sim c_1x^{-\alpha_    1}$ as $x\rightarrow\infty$ for $c_1>0$, $\alpha_1>1$. Then,
  $$\frac{1}{x\log x}\log\Pro\big(\cL(t)>x\big)\rightarrow-\alpha_1$$ as $x\rightarrow\infty.$
  \item[ii.)] Suppose that $\xi_0$ is such that $\Pro(\xi_0<-x)\sim c_2 x^{-\alpha_2}$ as $x\rightarrow\infty$ for $c_2>0$, $\alpha_2>1$. Then,
  $$\frac{1}{x\log x}\log\Pro\big(\cE(t)>x\big)\rightarrow-\alpha_2$$ as $x\rightarrow\infty.$
  \item[iii.)]Suppose that $\xi$ has a bounded density for which there exist positive constants $c_1,c_2,\alpha_1,\alpha_2$ such that
    \begin{align*}
      &f(x)  \sim c_1\,x^{-\alpha_1-1}, \\
      &f(-x)  \sim c_2\,x^{-\alpha_2-1}
    \end{align*}
    as $x\rightarrow\infty$. Then,
    $$\frac{1}{x\log x}\log \Pro\big(N(t)>x\big)\rightarrow -\min(\alpha_1+1,\alpha_2+1)$$ as $x\rightarrow\infty$.
\end{itemize}\label{Th: LD_L(t)_E(t)}
\end{theorem}
\hfill \break
\begin{proof}

      For part i.) we recall that $\cL(t)\overset{D}=\cL(0)$. Furthermore, $\Pro(\xi_{j+1}>j)\leq \Pro(-j+\xi_{-j}+U>0)\leq \Pro(\xi_j      >j-1)$, so that
  \begin{equation}\Pro\big(\sum_{j=1}^\infty I_j>x\big)\leq \Pro(\cL(t)>x)\leq \Pro\big(\sum_{j=0}^\infty I_j>x\big),\label{Eq: BoundsL(t)_Notes3.4}
  \end{equation}
where the $I_j$'s are independent Bernoulli r.v.'s in which $I_j =I(\xi_j>j-1)$. Theorem~\ref{Th: Bernoulli Tail} can then be applied to the extreme members of (\ref{Eq: BoundsL(t)_Notes3.4}), yielding i.). Part ii.) follows similarly. As for iii.), suppose that $\alpha_1\leq \alpha_2$ and set $I_j=I(j+\xi_j+U\in(0,t])$.  Fix an integer $d\geq 1$ and observe that

\begin{equation}
\begin{aligned}
\Pro\big(\sum_{j\leq 0}I_j>x\big)&\leq\Pro\big(N(t)>x\big)\\
&\leq \Pro\big(\sum_{j\leq 0}I_j>x\big)+\Pro\big(\sum_{j\leq 0}I_j\leq x,\sum_{j> 0}I_j>x\big)\\
&\leq\Pro\big(\sum_{j\leq 0}I_j>x\big)+\sum_{k=0}^{d-1}\Pro\Big(\sum_{j\leq 0}I_j\in x\, [k/d,(k+1)/d),\sum_{j> 0}I_j> x\,(1-(k+1)/d)\Big)\\
&\leq (d+1)\max_{0\leq k\leq d}\Pro\big(\sum_{j\leq 0}I_j\geq x\,k/d\big)\, \Pro\big(\sum_{j> 0}I_j\geq x\,(1-(k+1)/d)^+\big).
\end{aligned}\label{Eq: Bounds_N(t)_Notes3.5}
\end{equation}
Recalling our bounded density assumption, the Bounded Convergence Theorem implies that $\Pro(I_{-j}=1)\sim c_1\,t\,j^{-\alpha_1-1}$ and $\Pro(I_j=1)\sim c_2\,t\,j^{-\alpha_2-1}$ as $j\rightarrow\infty.$

Arguing as for i.) and ii.), we find that
\begin{equation}\label{Eq: LD1SumBern_Notes3.6}
\frac{1}{x\log x}\log\Pro\big(\sum_{j\leq 0}I_j>x\big)\rightarrow-(\alpha_1+1)
\end{equation}
and
\begin{equation}\label{Eq: LD2SumBern_Notes3.7}
\frac{1}{x\log x}\log\Pro\big(\sum_{j> 0}I_j>x\big)\rightarrow-(\alpha_2+1)
\end{equation}
as $x\rightarrow\infty$. Utilizing (\ref{Eq: LD1SumBern_Notes3.6}) and (\ref{Eq: LD2SumBern_Notes3.7}), we observe that, for $0\leq k\leq d-1$, the term
\begin{equation}\begin{aligned}
\frac{1}{x\log x}\log&\big[\Pro\big(\sum_{j\leq 0}I_j\geq x\,k/d\big)\, \Pro\big(\sum_{j> 0}I_j\geq x\,(1-(k+1)/d)\big]\\
&\rightarrow-\frac{(\alpha_1+1)\, k}{d}-\frac{(\alpha_2+1)\, (d-(k+1))}{d}\leq -(\alpha_1+1)\frac{d-1}{d}
\end{aligned}\end{equation}
as $x\rightarrow\infty$.
Hence, letting $x\rightarrow\infty$ on the extreme terms of (\ref{Eq: Bounds_N(t)_Notes3.5}), followed by sending $d\rightarrow\infty$ yields iii.). A symmetric argument works for $\alpha_1>\alpha_2$.


\end{proof}\\
The next result shows that the large deviations tail exponent of $N(t)$ is not inherited by its equilibrium limit. In other words, one cannot interchange $x\rightarrow\infty$ in the  large deviations limit with $t\rightarrow\infty$ in time.\\

For our next result, we fix $s\in[0,1]$ and recall Theorem~\ref{Th: N(n)-n} and the quantities ${\cal E}'(s)$ and ${\cal E}'(s)$ defined there.

\begin{prop}\label{prop:limit_interchange} Suppose that $\xi_0$ is such that $\Pro(\xi_0>x)\sim c_1\,x^{-\alpha_1}$ and $\Pro(\xi_0<-x)\sim c_2\,x^{-\alpha_2}$ as $x\rightarrow\infty$ for $c_1,c_2>0$ and $\alpha_1,\alpha_2>1$. Then, $$\frac{1}{x\log x}\log\Pro\Big(\big(\cE'(s)-\cL'(s)\big)-\big(\cE(0)-\cL(0)\big)>x\Big)\rightarrow-\min(\alpha_1,\alpha_2)$$ as $x\rightarrow\infty.$
\end{prop}
\hfill \break

\begin{proof}
  Utilizing Corollary~\ref{Cor: LD_Diff_Bern} and arguing as in the proof of Theorem~\ref{Th: LD_L(t)_E(t)}, we find that $$\frac{1}{x\log x}\log \Pro\big(\cE'(s)-\cL'(s)>x)\big)\rightarrow-\alpha_2,$$ and $$\frac{1}{x\log x}\log \Pro\big(\cL(0)-\cE(0)>x)\big)\rightarrow-\alpha_1$$ as $x\rightarrow\infty.$ We can now use the same upper bound argument as in (\ref{Eq: Bounds_N(t)_Notes3.5}) to conclude that
  $$\overline{\lim}_{x\rightarrow\infty}\,\frac{1}{x\log x}\log\Pro\Big(\big(\cE'(s)-\cL'(s)\big)-\big(\cE(0)-\cL(0)\big)>x\Big)\leq -\min(\alpha_1,\alpha_2).$$ For the lower bound, suppose that $\alpha_2\leq\alpha_1.$ We find that
  \begin{align*}
    \Pro\Big( &\big (\cE'(s)-\cL'(s)\big)-\big(\cE(0)-\cL(0)\big)>x\Big)  \\
    ~&\geq  \Pro\big( \cE'(s)-\cL'(s)>x+d\big)\,\Pro\big (\cL(0)-\cE(0)\geq-d\big),
  \end{align*} so that
  $$\underline{\lim}_{x\rightarrow\infty}\,\frac{1}{x\log x}\log\Pro\Big(\big(\cE'(s)-\cL'(s)\big)-\big(\cE(0)-\cL(0)\big)>x\Big)\geq -\alpha_2.$$
  A symmetric argument holds for $\alpha_1<\alpha_2.$
\end{proof}
\\

Because of their intrinsic interest and their importance for scheduled queues, we now provide exact tail asymptotics for Bernoulli sums.\\
\begin{theorem}~
\begin{itemize}
  \item [i.)] Suppose that there exists $c>0$ and $\alpha>1$ such that $p_n=c \,n^{-\alpha}(1+O(1/n))$ as $n\rightarrow\infty$. Then, $$\Pro(\sum_{j\geq 0}I_j\geq n)\sim \frac{1}{\sqrt{2\pi\eta_*}}\,r_*^{-n}\,n^{-\alpha n-1/2}\exp\big(\psi(r_* n^\alpha)\big)$$ as $n\rightarrow\infty,$ where $r_*$ satisfies $$\int_0^\infty \frac{c\,r_*}{c\,r_*+x^\alpha} dx =1$$ and $$\eta_*=\int_0^\infty \frac{c\,r_*\,x^\alpha}{(c\,r_*+x^\alpha)^2} dx.$$
    \item [ii.)] Suppose that $p_n=c(w+n)^{-\alpha}$ for $n\geq 0$, where $c, w>0$ and $\alpha>1$. Then,
    \begin{align*}
      \Pro(\sum_{j\geq 0}I_j\geq n)\sim &\Big(\frac{1}{\sqrt{2\pi}}\Big)^{\alpha+1}\,\frac{\Gamma(w)^\alpha}{\sqrt{\eta_*}}\,(c\,r_*)^{\frac{1}{2}-w}\,n^{-\alpha n+\frac{1}{2}\,(\alpha-1)-w\alpha}\,e^{\gamma n}
    \end{align*} as $n\rightarrow\infty,$ where $\gamma =\int_0^1\log\big(1+\frac{1}{c\,r_*}x^\alpha\big)dx+\int_1^\infty\log\big(1+c\,r_*\,x^{-\alpha}\big)dx +\alpha+\log(c)$.
\end{itemize}\label{Th: Exact_LD_Bernoulli}
\end{theorem}
\hfill \break
\begin{proof}
  We start from the change-of-measure formula (\ref{Eq: Change_measure}), with the specific choice $\theta_*(n)=\log r_*+\alpha\log n$ (so that $\exp(\theta_*(n))=r_* n^\alpha).$ Then, \begin{align*}\Pro(Z\geq n)=\exp\big(&-\theta_*(n)\,n+\psi(r_* n^\alpha)\big)\cdot\widetilde{\e}_nI(Z\geq n)\exp\big(-\theta_*(n)(Z-n)\big).\end{align*}
  We wish now to apply the local central limit theorem (CLT) to $Z$ under $\widetilde{P}_n.$ Recall (\ref{eq: Etilde_Z}) and note that 
  \begin{align*}
  \widetilde{\e}_nZ&=\sum_{j\geq 0} \frac{p_j\,r_*n^\alpha}{p_j\,(r_*\, n^\alpha-1)+1}\\
  &=\sum_{j\geq 0} \frac{(n/j)^\alpha c\,r_*}{(n/j)^\alpha\,c\, r_*-p_j+1}\,(1+O(1/j))\\
  &=\sum_{j\geq 0} \frac{ c\,r_*}{c\, r_* +(j/n)^\alpha- (j/n)^\alpha\,p_j}\,(1+O(1/j))\\
  &=\sum_{j\geq 0} \frac{ c\,r_*}{c\, r_* +(j/n)^\alpha + O(n^{-\alpha}) \,(1+O(1/j))}\\
  &=\sum_{j\geq 0} \frac{ c\,r_*}{c\, r_* +(j/n)^\alpha}(1+ O(1/j))\\
  &= \sum_{j\geq 0} \frac{c\,r_*}{c\, r_*+(j/n)^\alpha}+\sum_{0\leq j\leq k_n}O(1/j) +O(k_n^{-1})\sum_{j>k_n}\frac{c\,r_*}{c\, r_*+(j/n)^\alpha}\\
  &=n\,\sum_{j\geq 0}\frac{1}{n}\,v(j/n)+O(\log k_n)+O(n\,k_n^{-1})\,\sum_{j>k_n}\frac{1}{n}\,v(j/n),
  \end{align*}
where $v(x)=c\,r_*\,(c\,r_*+x^\alpha)^{-1}$ and $k_n$ is selected so that $k_n/n^{2/3}\rightarrow 1$ as $n\rightarrow\infty.$ But
\begin{align*}
n\sum_{j\geq 0}\frac{1}{n}\,v(j/n)=n\int_0^\infty v(x)dx+n\sum_{j\geq 0}\int_{j/n}^{(j+1)/n}[v(j/n)-v(x)]dx
\end{align*}
The defining equation for $r_*$ implies that $\int_0^\infty v(x)dx=1.$ 
Set $$\omega_n(x)=\int_{j/n}^x[v(j/n)-v(y)]dy.$$ Since $\omega_n$ is twice differentiable with $w'_n(x)=v(j/n)-v(x)$ and $w''_n(x)=-v'(x)$, there exists $x_{j,n}\in [j/n,(j+1)/n]$ such that $$\omega_n\big((j+1)/n\big)=\omega_n(j/n)+1/n\cdot \omega_n'(j/n)+1/n^2\cdot \omega_n''(x_{j,n})/2$$ so that $$\int_{j/n}^{(j+1)/n}[v(j/n)-v(y)]dy=-v'(x_{j,n})\cdot\frac{1}{2\,n^2}$$ and hence
$$n\sum_{j\geq 0}\int_{j/n}^{(j+1)/n}[v(j/n)-v(x)]dx=-1/2\,\sum_{j\geq 0}v'(x_{j,n})\,\frac{1}{n}.$$ The latter sum is a Riemann sum approximation to the integral of $-\frac{1}{2}v'(\cdot)$ over $[0,\infty).$ Consequently,
\begin{equation}n\,\sum_{j\geq 0}\int_{j/n}^{(j+1)/n}[v(j/n)-v(x)]dx\rightarrow-\frac{1}{2}\int_0^\infty v'(x)dx=1/2.\label{Eq: Riemann_Sum_Exact_Tail_Notes3.5a}\end{equation} Similarly, $$\sum_{j>k_n}\frac{1}{n}v(j/n)-\int_{k_n/n}^\infty v(x) dx\rightarrow 0$$ as $n\rightarrow\infty.$ It follows that
\begin{equation}\label{Eq: E_tildeZ_Notes3.6} \widetilde{\e}_nZ=n+O(n^{1/3})\end{equation} as $n\rightarrow\infty.$ Also, as noted in the proof of Theorem~\ref{Th: Bernoulli Tail},
\begin{equation}\label{Eq: Var_tildeZ_Notes3.7} \frac{1}{n}\,\tvar_nZ\rightarrow \eta_*\end{equation} as $n\rightarrow\infty.$

We are now ready to apply the local CLT due to Davis and McDonald (1995). We first write $Z=\sum_{j=0}^{b_n}I_j+Y_n$, where $b_n\rightarrow\infty$ fast enough that $\widetilde{\e}_nY_n^2/n\rightarrow 0$ as $n\rightarrow\infty.$ It is easily verified that the conditions of the Lindeberg-Feller CLT apply to $(Z-\widetilde{\e}_nZ)/(\tvar_nZ)^{1/2}$; see p. 205 of \cite{Chung74}. Furthermore,  by recalling that, $\tvar_nI_j=\widetilde{\Pro}_n(I_j=1)\,\widetilde{\Pro}_n(I_j=0)$, we conclude that the sequence $Q_n=\sum_{j\leq b_n}\min\big(\widetilde{\Pro}_n(I_j=0),\widetilde{\Pro}_n(I_j=1)\big)$ that appears in the hypotheses of Theorem~1.2 of \cite{DavisMcdon95} can be lower bounded by $\sum_{j\leq b_n}\tvar_nI_j\sim n\,\eta_*$ as $n\rightarrow\infty.$ Consequently, Theorem~1.2 asserts that
$$\widetilde{\Pro}_n(Z=k)=\phi\big(\frac{k-\widetilde{\e}_nZ}{\sqrt{n\,\eta_*}}\big)\,\frac{1}{\sqrt{n\,\eta_*}}\,(1+o(1))$$ uniformly in $k$ as $n\rightarrow\infty$, where $\phi(\cdot)$ is the density of a $\mathcal{N}(0,1)$ r.v. Hence, in view of (\ref{Eq: E_tildeZ_Notes3.6}) and (\ref{Eq: Var_tildeZ_Notes3.7}),
$$\widetilde{\Pro}_n(Z=n+k)=\phi\big(\frac{k}{\sqrt{n\,\eta_*}}\big)\,\frac{1}{\sqrt{n\,\eta_*}}\,(1+o(1))$$ as $n\rightarrow\infty,$ so that
\begin{align*}
&\widetilde{\e}_nI(Z\geq n)\exp\big(-\theta_*(n)\,(Z-n)\big)\\
&=\sum_{k\geq0}\widetilde{\Pro}_n(Z=n+k)\,\exp\big(-\theta_*(n)\,k\big)\\
&\sim \frac{1}{\sqrt{2\,\pi\,n\,\eta_*}}
\end{align*} as $n\rightarrow\infty$, proving part i.).

Part ii.) is a special case of i.). All that is needed is the development of an asymptotic for $\psi(r_*\,n^\alpha)$, to the order of $o(1).$ Denoting (as usual) the gamma function by $\Gamma(\cdot)$, we write

\begin{equation}\label{Eq: Psi(r*nalpha)_Notes3.8}
\begin{aligned}
\psi(r_*\,n^\alpha)&=\sum_{j=0}^{n-1}\log\big(c\,(j+w)^{-\alpha}\,(r_*\,n^{\alpha}-1)+1\big)\\
&~~~~~~+\sum_{j\geq n}\log\big(c\,(j+w)^{-\alpha}\,(r_*\,n^{\alpha}-1)+1\big)\\
&=\log\big(\prod_{j=0}^{n-1}\big(\frac{n}{j+w}\big)^\alpha(c\,r_*)^n\big)+\sum_{j=0}^{n-1}\log\big(1+\big(\frac{j+w}{n}\big)^\alpha\frac{1}{c\,r_*}-\frac{1}{c\,r_*\,n^\alpha}\big)\\
&~~~~~~+\sum_{j\geq n}\log\big(1+c\,r_*\big(\frac{j+w}{n}\big)^{-\alpha}(1-\frac{1}{r_*\,n^\alpha})\big)\\
&=\log\left(\left(\frac{n^n\,\Gamma(w)}{\Gamma(w+n)}\right)^\alpha(c\,r_*)^n\right)+n\int_{w/n}^{1+w/n}\log\big(1+\frac{1}{c\,r_*}\,x^\alpha-\frac{1}{c\,r_*n^\alpha}\big)dx\\
&~~~~~~+n\int_{1+w/n}^\infty\log\big(1+c\,r_*\,x^{-\alpha}\,(1-\frac{1}{r_*\,n^\alpha})\big)dx\\
&~~~~~~+n\,\sum_{j=0}^{n-1}\int_{j/n}^{(j+1)/n}[h_n(j/n)-h_n(x)]dx+n\,\sum_{j\geq n}\int_{j/n}^{(j+1)/n}[\tilde{h}_n(j/n)-\tilde{h}_n(x)]dx,
\end{aligned}
\end{equation}
where $$h_n(x)=\log\big(1+(x+\frac{w}{n})^\alpha\,\frac{1}{c\,r_*}-\frac{1}{c\,r_*n^\alpha}\big),$$
$$\tilde{h}_n(x)=\log\big(1+c\,r_*\,(x+\frac{w}{n})^{-\alpha}(1-\frac{1}{r_*\,n^\alpha})\big).$$

 The first term in the third equality is due to the property of the Gamma function whereby  for any $z>0$, $\Gamma(z+1)=z\Gamma(z)$.
Set $h(x)=\log\big( 1+\frac{1}{c\,r_*}\,x^\alpha\big),$ $\tilde{h}(x)=\log\big(1+c\,r_*\,x^{-\alpha}\big).$ Arguing as in (\ref{Eq: Riemann_Sum_Exact_Tail_Notes3.5a}), the sum of the last two terms converges to
\begin{equation}\label{Eq: logcr*_Notes3.9}\begin{aligned}
&-1/2\,\int_0^1h'(x)dx-1/2\,\int_1^\infty\tilde{h}'(x)dx\\
&=-1/2\,(h(1)-h(0))-1/2\,(\tilde{h}(\infty)-\tilde{h}(1))\\
&=-1/2\,(\log(1+\frac{1}{c\,r_*})-\log(1+c\,r_*))=1/2\,\log(c\,r_*).
\end{aligned}\end{equation}
Also,
\begin{equation}\label{Eq: log(1+1/cr*xalpha)_Notes3.10}\begin{aligned}
&n\,\int_{w/n}^{1+w/n}\log\big(1+\frac{1}{c\,r_*}\,x^\alpha-\frac{1}{c\,r_*n^\alpha}\big)dx\\
&=n\,\int_{w/n}^{1+w/n}\log\big(1+\frac{1}{c\,r_*}\,x^\alpha\big)dx+O(n^{1-\alpha})\\
&=n\,\int_{0}^{1}\log\big(1+\frac{1}{c\,r_*}\,x^\alpha\big)dx+w\,h(1)-w\,h(0)+o(1)
\end{aligned}
\end{equation}
as $n\rightarrow\infty.$ Similarly,
\begin{equation}\label{Eq: log(1+1/cr*xalpha)_Notes3.10}\begin{aligned}
&n\,\int_{1+w/n}^{\infty}\log\big(1+c\,r_*\,x^{-\alpha}\,(1-\frac{1}{n^\alpha})\big)dx\\
&=n\,\int_{1}^{\infty}\log\big(1+c\,r_*\,x^{-\alpha}\big)dx-w\,\tilde{h}(1).
\end{aligned}
\end{equation}
Finally, we use the asymptotic
$$\Gamma(w+n)\sim\sqrt{2\,\pi\,n}\left(\frac{w+n-1}{e}\right)^{n+w-1}$$ as $n\rightarrow\infty$ (see p. 63 of \cite{Feller71}), to conclude that
\begin{equation}\label{Eq: Gamma equiv_Notes3.12}\left(\frac{n^n\Gamma(w)}{\Gamma(w+n)}\right)^\alpha\sim\left(\frac{\Gamma(w)}{\sqrt{2\,\pi}}\right)^\alpha\,n^{(-w+1/2)\,\alpha}\,e^{\alpha\,n}\end{equation}
as $n\rightarrow\infty.$ Combining (\ref{Eq: Psi(r*nalpha)_Notes3.8}) through (\ref{Eq: Gamma equiv_Notes3.12}) yields part ii.).
\end{proof}\\

With Theorem~\ref{Th: Exact_LD_Bernoulli} at our disposal, we can now derive exact asymptotics for the r.v.'s $\cE(t)$ and $\cL(t)$. For example, in view of the fact that the proof of part $ii.)$ holds uniformly in $w$,
\begin{align*}
\mathbb{P}(\cE(0)\geq n)\sim & \Big(\frac{1}{\sqrt{2\pi}}\Big)^{\alpha+1}\,\sqrt{\frac{c\,r_*}{\eta_*}}\,\mathbb{E}\left[\frac{\Gamma(U)^\alpha}{(c\,r_*\,n^{\alpha})^{U}}\right]\,n^{-\alpha n+\frac{1}{2}\,(\alpha-1)}{\color{blue}\,e^{\gamma n}}
\end{align*} as $n\rightarrow\infty$, provided that $\mathbb{P}(\xi_0\leq -x)=c\,x^{-\alpha}$ for $x>1$ and where $\gamma$ is the same constant defined in Theorem~\ref{Th: Exact_LD_Bernoulli}

\section{Behavior of the $S/D/1$ Workload Process under Critical Loading}\label{sec: Critical Loading SD1}
In this section, we consider a queue that is fed by a scheduled traffic in which each customer's service time requirement is of unit duration, and in which the server has the capacity to process work at unit rate. Under these assumptions, the rate at which work arrives per unit time equals the service capacity of the system, so that the queue is subject to critical loading.

Note that the total work to arrive in $(0,t]$ is given by $N(t).$ Let $W(t)$ be the workload in the system at time $t$ (i.e. $t+W(t)$ is the first time subsequent to $t$ at which the system would empty if no additional work were to arrive after $t$.) If $W(0)=0,$ then $$W(t)=\max_{0\leq s\leq t}\Big[(N(t)-t)-(N(s)-s)\big].$$
Our goal is to analyze the behavior of $W(t)$ for $t$ large.

If $\mathbb{E}|\xi_0|<\infty$, Proposition~\ref{Prop: N(t)-t} applies  so that
\begin{equation}\label{Eq: Workload_Notes4.1}
W(t)=\max_{0\leq s\leq t}\Big[\cE(t)-\cE(s)-\cL(t)+\cL(s)\big]+O_p(1),
\end{equation}
where $O_p(1)$ is a term that is stochastically bounded in $t$. Because of the stationarity of $\big(\big(\cE(t),\cL(t)\big):t\in \mathbb{R}\big),$ the first term in (\ref{Eq: Workload_Notes4.1}) has the same distribution as
\begin{equation}\label{Eq: M(t)_Notes4.2}
\begin{aligned}
  M(t) & = \max_{0\leq s\leq t}\Big[\cE(0)-\cE(-(t-s))-\cL(0)+\cL(-(t-s))\Big],\\
  & =\max_{0 \leq r\leq t}\Big[\cE^*(r)-\cL^*(r)\Big]-\cE^*(0)+\cL^*(0),
\end{aligned}
\end{equation}
Note that $\cE^*(\cdot)$ is the ``early customer" process for the time-reversed system in which the perturbations are given by $(-\xi_{-j}:j\in\mathbb{Z}),$ and $\cL^*(\cdot)$ is the corresponding ``late customer" process. As a result, $\cE^*(r)\overset{d}=\cL(-r)$ and $\cL^*(r)\overset{d}=\cE(-r).$ As a matter of fact these  equalities hold pathwise except for the fact that  $\cE^*$ and $\cL^*$ are generated following the uniform distribution, $1-U$, instead of $U$. Given that $\cE^*$ and $\cL^*$ are non-negative processes for which $\cE^*(r)$ is independent of $\cL^*(r)$ for $r\in\mathbb{R}$, it is evident that the growth of $M(t)$ will be determined by $\cE^*$ and that the left tail of $-\xi_0$ (or right tail of $\xi_0$) governs the large time behavior of $M(\cdot)$ (and hence $W(\cdot))$. The dominance of the right tail of $\xi_0$ over the left tail is perhaps explained by the fact that the left tail induces the arrival of ``early customers" from the future evolution of the queue. More such early arrivals in an interval mean fewer potential customers available from which to stimulate a future burst of arrivals, so that the left tail has less influence over ``growing" $M(\cdot)$ over time.


\begin{theorem}\label{Th: 5_Workload_Crit_Load}
  Suppose that $\mathbb{E}\xi_0^-<\infty$ and that there exists constant $c>0$ and $\alpha>1$ for which $\mathbb{P}(\xi_0>x)\sim c\,x^{-\alpha}$ as $x\rightarrow\infty.$ Then, $$\frac{W(t)}{\log t/\log\log t}\Rightarrow 1/\alpha$$ as $t\rightarrow\infty.$
\end{theorem}
\hfill \break
\begin{proof}
Clearly,
\begin{align*}
  \max_{0\leq r\leq t} \big[\cE^*(r)-\cL^*(r)\big]& \leq\max_{0\leq r \leq t} \cE^*(r)\\
  & \leq \max_{1\leq n\leq \lfloor t\rfloor+1}\max_{0\leq s< 1}\cE^*(n-s).
  \end{align*}
For $0\leq s<1$ and $n\geq 1$,
\begin{align}
  \cE^*(n-s)& =\sum_{j+U\leq-n+s}I(j+U+\xi_j>-n+s)\nonumber\\
  & \leq \sum_{j+U\leq-n+s}I(j+U+\xi_j>-n)\nonumber\\
  &=\cE^*(n)+\sum_{-n<j+U\leq-n+s}I(j+U+\xi_j>-n)\nonumber\\
  &\leq \cE^*(n)+\sum_{-n<j+U\leq-n+s}1\nonumber\\
  &= \cE^*(n)+1\label{Eq: Ineq_E(n-s)_Notes4.2a}
  \end{align}
So,
$$\max_{0\leq r\leq t} \big[\cE^*(r)-\cL^*(r)\big]\leq 1+\max_{1\leq n\leq \lfloor t\rfloor+1}\cE^*(n).$$
Hence, for $\varepsilon>0$ and $t$ sufficiently large,
\begin{align}
&\mathbb{P}\Big(  \max_{0\leq r\leq t} \big[\cE^*(r)-\cL^*(r)\big]>\frac{1+3\,\varepsilon}{\alpha}\frac{\log t}{\log\log t}\Big)\nonumber\\
\leq &~\mathbb{P}\Big(  \max_{0\leq n\leq \lfloor t\rfloor+1} \cE^*(n)>\frac{1+2\,\varepsilon}{\alpha}\frac{\log t}{\log\log t}\Big)\nonumber\\
\leq &~\sum_{n=1}^{\lfloor t\rfloor+1}\mathbb{P}\Big(\cE^*(n)>\frac{1+2\,\varepsilon}{\alpha}\frac{\log t}{\log\log t}\Big)\nonumber\\
\leq &~(\lfloor t\rfloor+1)\,\mathbb{P}\Big(\cE^*(0)>\frac{1+2\,\varepsilon}{\alpha}\frac{\log t}{\log\log t}\Big)\nonumber\\
=&~(\lfloor t\rfloor+1)\,\exp\Big(\log\big(\mathbb{P}\big(\cL(0)>\frac{1+2\,\varepsilon}{\alpha}\frac{\log t}{\log\log t}\big)\big)\Big)\nonumber\\
\leq &~(\lfloor t\rfloor+1)\,\exp\big(-(1+\varepsilon)\,\log t\big)\nonumber\\
=&~O(t^{-\varepsilon})\rightarrow 0\label{Eq: O(t-eps)_Notes4.3}
  \end{align}
as $t\rightarrow \infty,$ where we used Theorem~\ref{Th: LD_L(t)_E(t)} for the final inequality.

To obtain the necessary lower bound, fix $\varepsilon\in(0,1/8\alpha)$ and note that for such $\varepsilon$, $1-2\,\varepsilon+\varepsilon^2<1-2\,\varepsilon-\varepsilon^2+\varepsilon/2\alpha\,(<1).$ Choose $\tau$ in the interval $(1-2\,\varepsilon+\varepsilon^2, 1-2\,\varepsilon-\varepsilon^2+\varepsilon/2\alpha)$. Put $b(t)=(1/\alpha)(\log t/\log\log t)$, $c(t)=(1-2\varepsilon)^2\,b(t)^2,$ and $k(t)=\lceil t^\tau\rceil$. As in the proof of Theorem~\ref{Th: Bernoulli Tail}, we find that for $\theta>0,$
\begin{align*}
\mathbb{P}(\cE(0)\geq n)&\leq \mathbb{P}\big(\sum_{j\geq -1}I(j+\xi_j\leq 0)\geq n\big)\\
&\leq \exp\Big(-\theta n+\sum_{j\geq -1}\log\big(\,\mathbb{P}(j+\xi_0\leq 0)(e^\theta-1)+1\big)\Big)\\
&\leq \exp\Big(-\theta n+\sum_{j\geq -1}\mathbb{P}(j+\xi_0\leq 0)(e^\theta-1)\Big)\\
&=\exp\Big(-\theta n+(e^\theta-1)\,\mathbb{E}\sum_{j= -1}^{\lceil -\xi_0\rceil}1\Big)\\
&\leq \exp\Big(-\theta n+(e^\theta-1)(\,\mathbb{E}\xi_0^-+3)\Big).
\end{align*}
By setting $\theta=\log n$, we conclude that $$\mathbb{P}\big(\cE(0)\geq n)\big)\leq \exp\big(-n \log n+O(n)\big)$$
so that
\begin{equation}\label{Eq: Inequality_E(0)Notes4.4}
  \mathbb{P}\big(\cE(0)\geq \varepsilon\,b(t)\big)\leq t^{-\varepsilon/2\alpha}
\end{equation}
for $t$ sufficiently large. In addition, an examination of the proof of Theorem~\ref{Th: Bernoulli Tail} shows that under the conditions stated there,
\begin{equation}\label{Eq: LogTail_limit_Notes4.5}
\frac{1}{z\log z} \log \mathbb{P}\big(\sum_{j=0}^{\lceil z^2\rceil} I_j>z\big)\rightarrow -\alpha
\end{equation}
as $z\rightarrow\infty.$

We now subdivide the interval $[-t,0]$ into $k(t)$ subintervals of equal length, and let $r_1,r_2,...,r_{k(t)}$ be the right endpoints of the $k(t)$ subintervals.

Then,
\begin{align}
&\mathbb{P}\Big(  \max_{0\leq r\leq t} \big[\cE^*(r)-\cL^*(r)\big]>(1-3\,\varepsilon)\, b(t)\Big)\nonumber\\
&\geq \mathbb{P}\Big(  \max_{1\leq i\leq k(t)} \big[\cL(r_i)-\cE(r_i)\big]>(1-3\,\varepsilon)\, b(t)\Big)\nonumber\\
&\geq \mathbb{P}\Big(  \max_{1\leq i\leq k(t)} \big[(\cL(r_i)-\cE(r_i))\,I\big(\cE(r_i)\leq \varepsilon\,b(t)\big)\big]>(1-3\,\varepsilon)\, b(t)\Big)\nonumber\\
&\geq\mathbb{P}\Big(  \max_{1\leq i\leq k(t)} \big[\cL(r_i)\,I\big(\cE(r_i)\leq \varepsilon\,b(t)\big)\big]>(1-2\,\varepsilon)\, b(t)\Big)\nonumber\\
&\geq \mathbb{P}\Big(\max_{1\leq i\leq k(t)}\cL(r_i)>(1-2\,\varepsilon)\,b(t)\Big)\nonumber \\
&~~~~~-\mathbb{P}\Big(\max_{1\leq i\leq k(t)}\cL(r_i)I\big(\cE(r_i)> \varepsilon\,b(t)\big)>(1-2\,\varepsilon)\, b(t)\Big)\nonumber\\
& \geq \mathbb{P}\Big(\max_{1\leq i\leq k(t)}\cL(r_i)>(1-2\,\varepsilon)\,b(t)\Big)\nonumber \\
&~~~~~-\sum_{i=1}^{k(t)}\mathbb{P}\Big(\cL(r_i)I\big(\cE(r_i)> \varepsilon\,b(t)\big)>(1-2\,\varepsilon)\, b(t)\Big)\nonumber\\
&\geq \mathbb{P}\Big(\max_{1\leq i\leq k(t)}\cL(r_i)>(1-2\,\varepsilon)\,b(t)\Big)\nonumber \\
&~~~~~-k(t)\,\mathbb{P}\Big(\cL(0)>(1-2\,\varepsilon)\, b(t),\cE(0)>\varepsilon\,b(t)\Big)\nonumber\\
&\geq \mathbb{P}\Big(\max_{1\leq i\leq k(t)}\sum_{r_i-c(t)\leq j\leq r_i-1} I(j+\xi_j>r_i)>(1-2\,\varepsilon)\,b(t)\Big)\nonumber \\
&~~~~~-k(t)\,\mathbb{P}\Big(\cL(0)>(1-2\,\varepsilon)\, b(t)\Big)\,\mathbb{P}\Big(\cE(0)>\varepsilon\,b(t)\Big)\nonumber\\
&= 1-\Big(1-\mathbb{P}\big(\sum_{-c(t)\leq j\leq -1} I(j+\xi_j>0)>(1-2\,\varepsilon)\,b(t)\big)\Big)^{k(t)}\label{Eq: term_Powerk(t)_Notes4.6}\\
&~~~~~-k(t)\,\mathbb{P}\Big(\cL(0)>(1-2\,\varepsilon)\, b(t)\Big)\,\mathbb{P}\Big(\cE(0)>\varepsilon\,b(t)\Big)\nonumber,
    \end{align}
where we used the independence of $\cL(0)$ and $\cE(0)$ and that of disjointly indexed indicator r.v.'s for both of the last two lines displayed above.  

Given (\ref{Eq: LogTail_limit_Notes4.5}), it follows that $$\mathbb{P}\Big(\sum_{-c(t)\leq j\leq -1}I(j+\xi_j>0)>(1-2\,\varepsilon)\,b(t)\Big)\geq t^{-(1-2\,\varepsilon)-\varepsilon^2}$$
for $t$ sufficiently large. In view of the choice of $\tau$, we conclude that
\begin{equation}\Big(1-\mathbb{P}\big(\sum_{-c(t)\leq j\leq -1} I(j+\xi_j>0)>(1-2\,\varepsilon)\,b(t)\big)\Big)^{k(t)}\rightarrow 0\label{Eq: Limitto_zero_Notes4.7}\end{equation}
as $t\rightarrow\infty$. On the other hand,
\begin{equation*}
  \mathbb{P}\Big(\cL(0)>(1-2\,\varepsilon)\, b(t)\Big)\leq t^{-(1-2\,\varepsilon)+\varepsilon^2}
\end{equation*}
for $t$ sufficiently large. Given (\ref{Eq: Inequality_E(0)Notes4.4}) and our choice of $\tau$, we find that
\begin{equation}\label{Eq: LimittozeroNotes4.8}
  k(t)\,\mathbb{P}\Big(\cL(0)>(1-2\,\varepsilon)\, b(t)\Big)\,\mathbb{P}\Big(\cE(0)>\varepsilon\,b(t)\Big)\rightarrow 0
\end{equation}
as $t\rightarrow\infty.$ Relations (\ref{Eq: O(t-eps)_Notes4.3}), (\ref{Eq: term_Powerk(t)_Notes4.6}), (\ref{Eq: Limitto_zero_Notes4.7}) and (\ref{Eq: LimittozeroNotes4.8}) prove the theorem.
\end{proof}\\

Theorem~\ref{Th: 5_Workload_Crit_Load} shows that the workload of the $S/D/1$ queue under critical loading increases very slowly (at $\log t/\big(\log \log t\big)$ rate), even in the presence of ``heavy tailed" perturbations. This is in sharp contrast to the $t^{1/2}$ increase in workload that occurs under critical loading for a $G/D/1$ queue, in which the arriving traffic is described by a renewal process with finite positive variance (see \cite{Glynn90}). This result makes clear the significant positive impact that scheduling can have upon queue performance.

\section{Behavior of the $S/D/1$ Workload Process in Heavy Traffic}\label{Sec:5_Heavy_Traffic}
We now turn to the analysis of the $S/D/1$ queue when the system has more service capacity than is needed. We assume, as in Section~\ref{sec: Critical Loading SD1}, that work is arriving at unit rate (on average) via deterministic service time requirements of unit size, but give the server a capacity to process work at the rate $1/\rho$ with $\rho<1$ (so that the queue's utilization factor is $\rho$). Let $W_\rho(\cdot)$ be the associated workload process. Then,
$$W_\rho(t)=\max_{0\leq s\leq t}\Big[\cE(t)-\cE(s)-\cL(t)+\cL(s)+a(t)-a(s)-\frac{1-\rho}{\rho}\,(t-s)\Big],$$
where $a(t)\overset{\Delta}=-(t-\lfloor t\rfloor)+I(U\leq t-\lfloor t\rfloor).$ As argued in Section~\ref{sec: Critical Loading SD1}, $W_\rho(t)\overset{D}=M_\rho(t),$ where
\begin{equation}\label{Eq: M_rho_Notes5.1}
M_\rho(t)=\max_{0\leq r\leq t}\Big[\cE^*(r)-\cL^*(r)-\frac{1-\rho}{\rho}\,r+a(0)-a(-r)\Big]+\cL^*(0)-\cE^*(0).
\end{equation}
Since $M_\rho(t)\nearrow M_\rho(\infty)$ a.s. as $t\rightarrow\infty$, it follows that $W_\rho(t)\Rightarrow W_\rho(\infty)$ as $t\rightarrow\infty$, where $W_\rho(\infty)\overset{D}= M_\rho(\infty)$. Our key result in this section describes the ``heavy traffic" behavior of $W_\rho(\infty)$ as $\rho\nearrow 1.$

\begin{theorem}\label{th: 6_HT_Result}
  Suppose that $\mathbb{E}\xi_0^-<\infty$ and that there exist constants $c>0$ and $\alpha>1$ for which $\mathbb{P}(\xi_0>x)\sim c\,x^{-\alpha}$ as $x\rightarrow\infty.$ Then,
  \begin{equation}\label{Eq: HT_thm_Notes5.2}
  \frac{\log \log \left(\frac{1}{1-\rho}\right)}{\log\left(\frac{1}{1-\rho}\right)}\,W_\rho(\infty)\,\Rightarrow\frac{1}{\alpha}
  \end{equation}
  as $\rho\nearrow 1.$
\end{theorem}
\hfill \break
\begin{proof}
Note that for $1/2<\rho<1$,
\begin{align*}
  &\max_{r\geq 0}\Big[\cE^*(r)-\cL^*(r) -\frac{1-\rho}{\rho}\,r\Big] \\
  \geq~&\max_{0\leq r\leq 1/(1-\rho)} \Big[\cE^*(r)-\cL^*(r) -\frac{1-\rho}{\rho}\,r\Big] \\
  \geq~&\max_{0\leq r\leq 1/(1-\rho)} \Big[\cE^*(r)-\cL^*(r) \Big]-2.
\end{align*}
Of course, Theorem~\ref{Th: 5_Workload_Crit_Load} establishes that
\begin{equation}\label{Eq: HT_thm_Notes5.3}
  \frac{\log \log \left(\frac{1}{1-\rho}\right)}{\log\left(\frac{1}{1-\rho}\right)}\max_{0\leq r\leq 1/(1-\rho)} \Big[\cE^*(r)-\cL^*(r) \Big]\Rightarrow\frac{1}{\alpha}
  \end{equation}
  as $\rho\nearrow 1$, proving the required lower bound for (\ref{Eq: HT_thm_Notes5.2}).

  To prove the upper bound, observe that
\begin{align}
&\max_{r\geq 0}\Big[\cE^*(r)-\cL^*(r) -\frac{(1-\rho}{\rho}\,r\Big]\nonumber\\
\leq~&\max_{0\leq r\leq \left(\frac{1}{1-\rho}\right)^{1+\varepsilon}} \Big[\cE^*(r)-\cL^*(r) \Big]+\max_{r\geq\left(\frac{1}{1-\rho}\right)^{1+\varepsilon}}\Big[\cE^*(r)-\frac{1-\rho}{\rho}\,r\Big]\label{Eq: UBD_pf_Th6_Notes5.4}.
\end{align}

Application of Theorem~\ref{Th: 5_Workload_Crit_Load} proves that
\begin{equation}\label{Eq: WeakLimit_prrofTh6_Notes5.5}
  \frac{\log \log \left(\frac{1}{1-\rho}\right)}{\log\left(\frac{1}{1-\rho}\right)}\max_{0\leq r\leq\left(\frac{1}{1-\rho}\right)^{1+\varepsilon}} \Big[\cE^*(r)-\cL^*(r) \Big]\Rightarrow\frac{1+\varepsilon}{\alpha}
\end{equation}
 as $\rho\nearrow 1$. On the other hand,
 \begin{align}
&\mathbb{P}\Big( \max_{r\geq\left(\frac{1}{1-\rho}\right)^{1+\varepsilon}}\Big[\cE^*(r)-\frac{1-\rho}{\rho}\,r\Big]\geq 1\Big)\nonumber\\
\leq~& \mathbb{P}\Big(\max_{n\geq 0}\Big[ \max_{0\leq s\leq 1}\cE^*\Big(\big(\frac{1}{1-\rho}\big)^{1+\varepsilon}+n+s\Big)-\big(\frac{1}{1-\rho}\big)^\varepsilon\,\frac{1}{\rho}-\frac{1-\rho}{\rho}\,n\Big]\geq 1\Big)\nonumber\\
\leq~& \mathbb{P}\left(\max_{n\geq 1}\Big[ \cE^*\Big(\big(\frac{1}{1-\rho}\big)^{{\color{blue}1+}\varepsilon}+n\Big)-\frac{1}{(1-\rho)^\varepsilon}-(1-\rho)\,n\Big]\geq 0\right),\label{Eq: Ineq_Proofth6_Notes5.6}
\end{align}
where we used (\ref{Eq: Ineq_E(n-s)_Notes4.2a}) for the last inequality. The quantity (\ref{Eq: Ineq_Proofth6_Notes5.6}) can, in turn, be upper bounded by
$$\sum_{n=0}^\infty\mathbb{P}\Big(\cE^*(0)\geq \frac{1}{(1-\rho)^\varepsilon}+(1-\rho)\,n\Big).$$
Theorem~\ref{Th: LD_L(t)_E(t)} proves that $$\mathbb{P}(\cE^*(0)\geq t) \leq \exp(-\alpha\,t)$$ for $t$ sufficiently large, and hence the above sum is dominated by
\begin{align*}
&\sum_{n=0}^\infty\exp\Big(-\alpha\,\big(\frac{1}{(1-\rho)^\varepsilon}+(1-\rho)\,n\big)\Big)\\
=~&\exp\Big(-\,\frac{\alpha}{(1-\rho)^\varepsilon}\Big)\,\Big(1-\exp\big(-\alpha\,(1-\rho)\big)\Big)^{-1}\\
\sim~& \exp\Big(-\,\frac{\alpha}{(1-\rho)^\varepsilon}\Big)\,\big(\frac{1}{\alpha\,(1-\rho)}\big)\rightarrow 0
\end{align*}
as $\rho\nearrow 1$. Relations (\ref{Eq: HT_thm_Notes5.3}), (\ref{Eq: UBD_pf_Th6_Notes5.4}), (\ref{Eq: WeakLimit_prrofTh6_Notes5.5}), and (\ref{Eq: Ineq_Proofth6_Notes5.6}) then prove the theorem, in view of the fact that $\varepsilon$ can be made arbitrarily small.
\end{proof}\\

This $S/D/1$ heavy traffic limit theorem should be contrasted against the analogous $G/D/1$ limit theorem, for which the steady-state r.v. $W_\rho(\infty)$ scales as $1/(1-\rho)$ as $\rho\nearrow 1$; (see \cite{Glynn90}). For the $G/D/1$ queue, time scales of order $1/(1-\rho)^2$ are needed in order that fluctuations of order $1/(1-\rho)$ are exhibited (when $W_\rho(0)=0)$ (see again \cite{Glynn90}). The proof of Theorem~\ref{th: 6_HT_Result} shows that the time scale needed for $W_\rho$ to reach equilibrium is of order $1/(1-\rho)$, so that the $S/D/1$ queue equilibrates more quickly than does the $G/D/1$ queue.
\section{Remarks on the S/G/1 Queue}\label{sec: 6_SG1}
This paper has focussed thus far on the $S/D/1$ queue. We now turn to a discussion of the $S/G/1$ queue, in which the service requirements $(V_i:i\geq 1)$ associated with the sequence of arriving customers is assumed to be i.i.d. and independent of $N$. In this setting the total work $\Lambda(t)$ to arrive in the interval $(0,t]$ is given by $$\Lambda(t)=\sum_{i=1}^{N(t)}V_i.$$ Our main objective here is to point out that if the $V_i$'s are random (i.e. non-degenerate), then the behavior of the $S/G/1$ queue closely resembles that of the corresponding $D/G/1$ queue in which $\xi_i\equiv0$ for $i\in \mathbb{Z}$.

Our first result shows that, in great generality, $\Lambda$ satisfies the same functional central limit theorem (FCLT) as does $\Lambda',$ where $$\Lambda'(t)=\sum_{i=1}^{\lfloor t \rfloor} V_i.$$\\

\begin{prop}\label{prop: 4_S/G/1}Suppose that $\mathbb{E}|\xi_0|<\infty$ and $\mathbb{E}V_1^p<\infty$ for $p>2$. Then
$$\frac{1}{\sqrt{t}}\,\sup_{0\leq s\leq t}|\Lambda(s)-\Lambda'(s)|\Rightarrow 0$$ as $t\rightarrow\infty.$
\end{prop}
\hfill \break

\begin{proof}
First, we recall that $\max_{1\leq i \leq n}|V_i|=o(n^{1/p})$ a.s. as $n\rightarrow\infty$; (see, p. 278 of \cite{Feller71}). Also,
\begin{align}
&t^{-1/2}\max_{0\leq s\leq t}\left|\sum_{i=1}^{N(s)}V_i-\sum_{i=1}^{\lfloor s \rfloor} V_i\right|\nonumber\\
&\leq t^{-1/2}\max_{0\leq i\leq N(t)+t}\left|V_i\right|\cdot\max_{0\leq s\leq t}|N(s)-s|\nonumber\\
&= \frac{\max_{0\leq i\leq 3\,t}\left|V_i\right|}{t^{1/p}}\cdot\frac{\max_{0\leq s\leq t}|N(s)-s|}{t^{1/2-1/p}}\nonumber\\
&=o(1)\cdot \frac{\max_{0\leq s\leq t}|N(s)-s|}{t^{1/2-1/p}}~~a.s.
\end{align}
as $t\rightarrow\infty,$  where the first equality is due to the fact $N(t)-t\leq t$ a.s. We further note that $\mathbb{E}\xi_0<\infty$ guarantees that $\mathbb{E}\exp\big(\theta\cL(n)\big)<\infty$ for $\theta>0$. In particular, $\sum_{n\geq 1}\mathbb{P}\big(\exp(\theta\cL(n))>n\big)<\infty$, and hence the Borel-Cantelli lemma insures that, $\overline{\lim}_{n\rightarrow\infty}\frac{\cL(n)}{\log n}\leq 1/\theta~~a.s.$ Hence,$$ \overline{\lim}_{n\rightarrow\infty}\frac{\max_{1\leq i\leq n}\cL(i)}{\log n}\leq 1/\theta~~a.s.$$

Since $\left|\cL(n+s)-\cL(n)\right|\leq 1$ for $0\leq s<1$ (see (\ref{Eq: Ineq_E(n-s)_Notes4.2a}) for a similar bound involving $\cE$),
$$ \overline{\lim}_{t\rightarrow\infty}\frac{\max_{1\leq s\leq t}\cL(s)}{\log t}\leq 1/\theta~~a.s.$$ for $\theta>0$. Similarly, $\overline{\lim}_{t\rightarrow\infty}\max_{1\leq s\leq t}\cE(s)/\log t\leq 1/\theta~~a.s.$ for $\theta>0$, because $\mathbb{E}\xi_0^-<\infty.$ Proposition~\ref{Prop: N(t)-t} then completes the proof.
\end{proof}\\

 As a consequence, the input to the $S/G/1$ queue satisfies the same FCLT as for the $D/G/1$. Hence, the heavy traffic theory for the $S/G/1$ with random service times is identical to that for the corresponding $D/G/1$ queue.  This fact is illustrated in the next result.

\begin{corollary}
  Under the assumptions of Proposition~\ref{prop: 4_S/G/1}
  $$n^{-1/2}\Big(\sum_{i=1}^{N(n\cdot)}V_i-n\,\mathbb{E}V_1\cdot\big)\Rightarrow\sqrt{\textnormal{\var} V_1}\,B(\cdot)$$
as $n\rightarrow\infty$, in $D(0,\infty)$, where $B(\cdot)$ is a standard Brownian motion, and $\Rightarrow$ corresponds here to weak convergence in $D(0,\infty)$ (see, \cite{Billingsley99} for the definition).
\end{corollary}

\bibliographystyle{ormsv080}
\bibliography{BiblioApP}
\end{document}